\documentclass{article}
\usepackage{mathtools,amsmath,amssymb,amsthm}
\usepackage[lite]{amsrefs}
\usepackage{microtype}

\theoremstyle{plain}
\newtheorem{theorem}[equation]{Theorem}
\newcommand*{\defeq}{\mathrel{\vcentcolon=}}
\newcommand*{\Cst}{\mathrm C^*}
\newcommand*{\K}{\mathrm K}
\newcommand*{\Comp}{\mathbb K}
\newcommand*{\Mat}{\mathbb M}

\title{Topological insulators and stable isomorphism versus
  isomorphism of vector bundles}
\author{Ralf Meyer\thanks{This work was supported by the Shota
  Rustaveli National Science Foundation of Georgia (SRNSFG) grant
  FR-23-779.}}

\begin{document}
\maketitle

\begin{abstract}
  This note gives an overview of the mathematical framework underlying
  topological insulators, highlighting the connection to $\K$-theory
  and vector bundles.  We see ``real'' and ``quaternionic'' vector
  bundles arise naturally in the presence of time-reversal symmetry.
  Our recent results about when stable isomorphism implies isomorphism
  are summarised, including some ongoing work for $G$-equivariant
  $\K$-theory for finite groups.  This clarifies when \(\K\)-theory
  completely distinguishes topological phases. 
\end{abstract}

\section{Introduction}

Topological insulators have become a central topic in condensed matter
physics because they exhibit remarkable transport properties: they are
insulating in the bulk, but they support robust conducting boundary
states, which are protected by topological invariants of the system.
The electronic properties of these materials are determined by global
topological invariants, making them resilient to perturbations.
Mathematically, they may be described in the language of operator
algebras, $\K$-theory, and vector bundles.

In this note, we briefly review the mathematical formalism that models
topological insulators, emphasizing connections to $\Cst$-algebras and
vector bundles.  We describe how time-reversal symmetry introduces
``real'' and ``quaternionic'' structures.  We explain the physical
motivation for the question whether stably isomorphic vector bundles
are isomorphic, and we discuss our recent results on this in joint
work with Malkhaz Bakuradze~\cite{Bakuradze-Meyer:Iso_stable_iso}.  We
conclude with some remarks on ongoing research in $G$-equivariant
$\K$-theory, which exhibits markedly different phenomena compared to
the non-equivariant case.

\section{Lattice models and the Hilbert space formalism}

We consider an electron moving in a $d$-dimensional crystal with $k$
internal degrees of freedom per lattice site.  The Hilbert space of
the system is \(\ell^{2}(\mathbb{Z}^{d},\mathbb{C}^{k})\).  The
dynamics on this space is governed by a bounded, self-adjoint
Hamiltonian~$H$.  We assume translation invariance, that is,
\(S_n H = H S_n\) for all \(n \in \mathbb{Z}^d\) for the translation
operator~$S_n$ defined by $(S_n f)(x) = f(x-n)$.  Under this
assumption, there are matrices $H_a \in \Mat_k(\mathbb{C})$ with
\begin{equation}
  \label{eq:Hamiltonian_coeff}
  (H g)(m) = \sum_{a \in \mathbb{Z}^d} (H_a g)(m-a)
\end{equation}
for all \(g\in \ell^2(\mathbb Z^d,\mathbb C^k)\).
We assume that~$H$ has finite range, that is, $H_a = 0$ for all but
finitely many~$a$.  A Hamiltonian~$H$ describes an insulator if and
only if it is invertible.

To describe a crystal with a boundary, we restrict~$H$ to a
half-space, that is, we consider the operator
\(\hat{H} \defeq I^* H I\) on the subspace
\(\ell^2(\mathbb{N} \times \mathbb{Z}^{d-1}, \mathbb{C}^k)\), where
\(I\colon \ell^2(\mathbb{N} \times \mathbb{Z}^{d-1}, \mathbb{C}^k)
\hookrightarrow \ell^2(\mathbb{Z}^{d}, \mathbb{C}^k)\) denotes the
inclusion map.  Even if~\(H\) is invertible, \(\hat{H}\) may fail to
be invertible.  This means physically that conducting states appear on
the boundary of a finite-size chunk of the material.  Next, we sketch
how the presence of such conducting states may be topologically
protected by a non-vanishing $\K$-theory index.

\section{C*-Algebraic formulation and the index map}

The $\Cst$-algebra generated by the finite-range, translation invariant
operators on \(\ell^{2}(\mathbb{Z}^{d},\mathbb{C}^{k})\) is
\[
\Cst(\mathbb{Z}^d) \otimes \Mat_k(\mathbb{C}),
\]
where $\Cst(\mathbb{Z}^d) \cong C(\mathbb{T}^d)$ is the group
$\Cst$-algebra of $\mathbb{Z}^d$.  This is the smallest $\Cst$-algebra
that contains all the Hamiltonians allowed above.  On the half-space,
one unitary generator is replaced by a unilateral shift on
\(\ell^2(\mathbb N)\).  So the allowable Hamiltonians for the
half-space generate the $\Cst$-algebra
\[
\mathcal{T} \otimes \Cst(\mathbb{Z}^{d-1}) \otimes \Mat_k(\mathbb{C}),
\]
where $\mathcal{T}$ is the Toeplitz $\Cst$-algebra
(see~\cite{Prodan-Schulz-Baldes:Bulk_boundary} for more details).
These $\Cst$-algebras fit in an extension
\begin{equation}
  \label{eq:extension}
  \Comp \otimes \Cst(\mathbb{Z}^{d-1}) \otimes \Mat_k(\mathbb{C}) 
  \rightarrowtail \mathcal{T} \otimes \Cst(\mathbb{Z}^{d-1}) \otimes \Mat_k(\mathbb{C}) 
  \twoheadrightarrow \Cst(\mathbb{Z}^d) \otimes \Mat_k(\mathbb{C}),
\end{equation}
where~$\Comp$ is the $\Cst$-algebra of compact operators on
$\ell^2(\mathbb N)$.  Roughly speaking, a half-space Hamiltonian in
$\mathcal{T} \otimes \Cst(\mathbb{Z}^{d-1}) \otimes \Mat_k(\mathbb{C})$
yields a bulk Hamiltonian in
$\Cst(\mathbb{Z}^d) \otimes \Mat_k(\mathbb{C})$ by identifying
$\mathbb{N} \times \mathbb{Z}^{d-1}$ with
$(\mathbb{N}-s) \times \mathbb{Z}^{d-1}$ for \(s\in\mathbb N\) and letting
$s\to\infty$.

Next we bring $\K$-theory and the index map into play.  The Hamiltonian
of an insulator is a self-adjoint invertible element of the
$\Cst$-algebra $\Cst(\mathbb{Z}^d) \otimes \Mat_k(\mathbb{C})$.  Functional
calculus allows us to deform it among self-adjoint invertible
operators to the self-adjoint involution
\(F = \operatorname{sign}(H)\).  This satisfies $F^2 = 1$ and
$F = F^*$.  It contains the same information as the associated
projection
\[
  p = \frac{1 + F}{2}
  =  \frac{1 + \operatorname{sign}(H)}{2}
  \in \Cst(\mathbb{Z}^d) \otimes \Mat_k(\mathbb{C}).
\]
The latter represents a class
$[H] \defeq [p]\in \K_0(\Cst(\mathbb{Z}^d))$.  The boundary map for
the extension~\eqref{eq:extension} maps this class to its index in
$\K_1(\Comp \otimes \Cst(\mathbb{Z}^{d-1}))$.  If the latter is
nonzero, then~$H$ cannot lift to an invertible operator in
$\mathcal{T} \otimes \Cst(\mathbb{Z}^{d-1}) \otimes
\Mat_k(\mathbb{C})$.  This means that any Hamiltonian on the
half-space that behaves like $H$ in the bulk is a conductor.  That is,
there are conducting states on the boundary, and these are forced to
exist by the nonvanishing index.  Since the index is homotopy
invariant, the existence of these boundary states is not affected by
small perturbations of the Hamiltonian.

\section{Vector bundles and the Bloch bundle}

Via Fourier transform, $\Cst(\mathbb{Z}^d)$ is identified with
$C(\mathbb{T}^d)$, so that the $\K_0$-class $[H]$ defines a class in
$\K^0(\mathbb{T}^d)$.  The associated vector bundle is called the
\emph{Bloch bundle}.  Its fibre at $z\in \mathbb{T}^d$ is the image of
the Fourier transform of~$p$ at $z$; the latter is a projection in
$\Mat_k(\mathbb{C})$.

It is physically interesting to know whether or not this vector bundle
is trivial because this is equivalent to the existence of
``exponentially localised Wannier functions'' (see
\cite{DeNittis-Lein:Wannier_zero_flux}*{Proposition~4.3}), which are a
tool used for computations in physics.  It is usually much easier to
decide whether the Bloch bundle is stably trivial, that is, becomes
trivial after adding a trivial bundle.  This means that its class in
reduced \(\K\)-theory vanishes.  If the reduced \(\K\)-theory is
torsion-free, this happens if and only if its Chern numbers vanish.
Physicists have long studied the Chern numbers of the Bloch bundle as
topological invariants related to conductivity phenomena.  As a
result, it is physically relevant to know whether the triviality of
the Bloch bundle follows from its stable triviality.  It is well known
that vector bundles of sufficiently high rank that are stably trivial
are automatically trivial; more generally, stably isomorphic bundles
of sufficiently high rank are isomorphic (see~
\cite{Husemoller:Fibre_bundles_3}).  Generalizations of these
classical results are needed to treat physical systems with certain
extra symmetries.

\section{``Real'' and ``Quaternionic'' bundles}

In quantum mechanics, a time-reversal symmetry is represented by an
anti-unitary operator that commutes with the Hamiltonian.  The square
of this anti-unitary operator is $\pm1$, where $+1$ occurs for bosons,
and~$-1$ for fermions.  In the Hilbert space
$\ell^2(\mathbb{Z}^d, \mathbb{C}^k)$, we assume that time-reversal
symmetry is given by applying a certain anti-unitary operator
on~$\mathbb{C}^k$ pointwise.  Depending on the sign, this has the
effect that the coefficients~$H_a$ in~\eqref{eq:Hamiltonian_coeff} now
belong to $\Mat_k(\mathbb{R})$ or $\Mat_k(\mathbb{H})$ instead of
$\Mat_k(\mathbb{C})$.  The same happens for the matrix coefficients of
the projection~$p$.  When we take the Fourier transform, however, this
does not correspond to the Bloch bundle over the torus being a real or
quaternionic vector bundle.  Instead, the Bloch bundle is a complex vector
bundle equipped with a conjugate-linear involution $\theta$ mapping
$E_z$ to $E_{\bar{z}}$ for all $z\in \mathbb{T}^d$.  If $\theta^2 = 1$
(bosons), this is a ``real'' vector bundle; if $\theta^2 = -1$, this
is a ``quaternionic'' vector bundle.

Such vector bundles may be consideded over a space~\(X\) with an
involution such as the map \(z\mapsto \bar{z}\) above.  If the
involution on~\(X\) is the identity map, then ``real'' and
``quaternionic'' vector bundle become equivalent to vector bundles
over the fields \(\mathbb R\) of real numbers and~\(\mathbb H\) of
quaternions, respectively.  We need the case, however, where the
involution is nontrivial.  While ``real'' vector bundles have been
known in index theory for a long time (see~\cite{Atiyah:K_Reality}),
they have not received so much attention.  In particular, it has not
been shown that ``real'' or ``quaternionic'' bundles of sufficiently
high rank are isomorphic once they are stably isomorphic.  The recent
article~\cite{Bakuradze-Meyer:Iso_stable_iso} fills this gap.  A
crucial step in the proof is showing that a bundle of sufficiently
high rank contains a trivial vector bundle of rank~$1$ as a direct
summand.  Both results combine into the statement that the
stabilisation map from bundles of rank~$k$ to bundles of rank~$k+1$
that adds a trivial bundle of rank~$1$ induces a bijection on
isomorphism classes for sufficiently large~$k$.  The following
theorems from~\cite{Bakuradze-Meyer:Iso_stable_iso} give the details
of these statements:

\begin{theorem}
  \label{the:main_real}
  Let \(d_1,d_0,k\in\mathbb{N}\).  Let~\(X\) be a
  \(\mathbb{Z}/2\)-CW-complex.  Assume that the free cells in~\(X\)
  have at most dimension~\(d_1\) and that the trivial cells have at
  most dimension~\(d_0\).  Let
  \[
    k_0 \defeq
    \max \left\{ d_0,
      \left\lceil \frac{d_1-1}{2} \right\rceil \right\},
    \qquad
    k_1 \defeq
    \max \left\{ d_0+1,
       \left\lceil \frac{d_1}{2} \right\rceil \right\}.
  \]
  \begin{enumerate}
  \item \label{the:main_real_1}%
    Let~\(E\) be a ``real'' vector bundle over~\(X\) of rank
    \(k\ge k_0\).  There is an isomorphism
    \(E\cong E_0 \oplus (X\times \mathbb{C}^{k-k_0})\) for some ``real''
    vector bundle~\(E_0\) over~$X$ and the trivial ``real'' vector
    bundle $X\times \mathbb{C}^{k-k_0}$ of rank~$k-k_0$.
  \item \label{the:main_real_2}%
    Let \(E_1\) and~\(E_2\) be two ``real'' vector bundles over~\(X\)
    of rank \(k\ge k_1\).  If \(E_1\) and~\(E_2\) are stably
    isomorphic, that is, $E_1 \oplus E_3 \cong E_2 \oplus E_3$ for
    some ``real'' vector bundle~$E_3$, then they are isomorphic.
  \end{enumerate}
\end{theorem}

\begin{theorem}
  \label{the:main_quaternion}
  Let \(d_1,d_0,k\in\mathbb{N}\).  Let~\(X\) be a
  \(\mathbb{Z}/2\)-CW-complex.  Assume that the free cells in~\(X\)
  have at most dimension~\(d_1\) and that the trivial cells have at
  most dimension~\(d_0\).  Let
  \[
    k_0 \defeq \max \left\{ \left\lceil\frac{d_0-3}{2}\right\rceil,
      \left\lceil \frac{d_1-1}{2} \right\rceil \right\},\qquad
    k_1 \defeq \max \left\{ \left\lceil\frac{d_0-2}{2} \right\rceil,
    \left\lceil \frac{d_1}{2} \right\rceil \right\}.
  \]
  \begin{enumerate}
  \item \label{the:main_quaternion_1}%
    Let~\(E\) be a ``quaternionic'' vector bundle over~\(X\) of rank
    \(k\ge k_0\).  There is an isomorphism
    \(E\cong E_0 \oplus \theta^{\oplus 2\lfloor (k-k_0)/2\rfloor}_X\)
    for some ``quaternionic'' vector bundle~\(E_0\) and the trivial
    ``quaternionic'' vector bundle
    $\theta^{\oplus 2\lfloor (k-k_0)/2\rfloor}_X$ of
    rank~$2\lfloor (k-k_0)/2\rfloor$.
  \item \label{the:main_quaternion_2}%
    Let \(E_1\) and~\(E_2\) be two ``quaternionic'' vector bundles
    over~\(X\) of rank \(k\ge k_1\).  If \(E_1\) and~\(E_2\) are
    stably isomorphic, that is, $E_1 \oplus E_3 \cong E_2 \oplus E_3$
    for some ``quaternionic'' vector bundle~$E_3$, then they are
    isomorphic.
  \end{enumerate}
\end{theorem}

\section{Equivariant $\K$-theory and ongoing work}

When the system has classical crystallographic symmetries, then the
Bloch bundle becomes an equivariant vector bundle for a certain group.
In $G$-equivariant $\K$-theory for, say, a finite group~$G$, new
phenomena may occur.  The main new issue is that there may be several
non-isomorphic trivial vector bundles, corresponding to inequivalent
representations of~$G$.  This is not the case for ``real'' and
``quaternionic'' bundles, although they may at first sight seem more
complicated than equivariant vector bundles because they involve the
group $\mathbb{Z}/2$ acting on the vector bundle by anti-unitary maps.

For example, consider a $\mathbb{Z}/2$-equivariant complex vector
bundle over the circle
$\mathbb{T}^1\defeq \{z\in\mathbb{C}\mid \lvert z \rvert = 1\}$ with
$\mathbb{Z}/2$ acting by $z \mapsto \bar{z}$.  Its fibres over~$\pm1$
are complex representations of the group~$\mathbb{Z}/2$, and both may
be arbitrary.  It may happen that at~\(+1\) we have the trivial
representation and at~\(-1\) the nontrivial sign representation of
some rank~$k$.  This vector bundle has arbitrarily high rank~$k$, but
has no trivial subbundles because the representations at $\pm1$ have
no common subrepresentation.  Analogues of the main theorems above
exist, but their assumptions use multiplicities of representations of
stabiliser subgroups instead of just ranks.

\section{Conclusion}

Topological insulators illustrate a rich interplay between condensed
matter physics and $\K$-theory. Their mathematical description via
$\Cst$-algebras, vector bundles, and index maps explains the
topological protection of boundary states and motivates questions on
when stable isomorphism implies isomorphism for ``real'',
``quaternionic'', and equivariant bundles.  We stated two theorems
from~\cite{Bakuradze-Meyer:Iso_stable_iso} about trivial direct
summands and stable isomorphism and isomorphism of ``real'',
``quaternionic'' vector bundles.


\bib{Atiyah:K_Reality}{article}{
  author={Atiyah, Michael F.},
  title={$K$\nobreakdash -theory and reality},
  journal={Quart. J. Math. Oxford Ser. (2)},
  volume={17},
  date={1966},
  pages={367--386},
  issn={0033-5606},
  review={\MR {0206940}},
  doi={10.1093/qmath/17.1.367},
}

\bib{Bakuradze-Meyer:Iso_stable_iso}{article}{
  author={Bakuradze, Malkhaz},
  author={Meyer, Ralf},
  title={Isomorphism and stable isomorphism in ``real'' and ``quaternionic'' K-theory},
  journal={New York J. Math.},
  volume={31},
  pages={690--700},
  issn={1076-9803},
  review={\MR {4895179}},
  eprint={https://nyjm.albany.edu/j/2025/31-25.html},
  date={2025},
}

\bib{Husemoller:Fibre_bundles_3}{book}{
  author={Husemoller, Dale},
  title={Fibre bundles},
  series={Graduate Texts in Mathematics},
  volume={20},
  edition={3},
  publisher={Springer-Verlag},
  place={New York},
  date={1994},
  pages={xx+353},
  isbn={0-387-94087-1},
  review={\MR {1249482}},
}

\bib{DeNittis-Lein:Wannier_zero_flux}{article}{
  author={De Nittis, Giuseppe},
  author={Lein, Max},
  title={Exponentially localized Wannier functions in periodic zero flux magnetic fields},
  journal={J. Math. Phys.},
  volume={52},
  date={2011},
  number={11},
  pages={112103, 32},
  issn={0022-2488},
  review={\MR {2906554}},
  doi={10.1063/1.3657344},
}

\bib{Prodan-Schulz-Baldes:Bulk_boundary}{book}{
  author={Prodan, Emil},
  author={Schulz-Baldes, Hermann},
  title={Bulk and boundary invariants for complex topological insulators},
  series={Mathematical Physics Studies},
  note={From $K$-theory to physics},
  publisher={Springer},
  date={2016},
  pages={xxii+204},
  isbn={978-3-319-29350-9},
  isbn={978-3-319-29351-6},
  review={\MR {3468838}},
  doi={10.1007/978-3-319-29351-6},
}



\begin{bibdiv}
  \begin{biblist} \bibselect{references}
  \end{biblist}
\end{bibdiv}
\end{document}